\documentclass[11pt]{gtart}
\usepackage{amsmath,amssymb}
\usepackage[british]{babel}
\usepackage[all]{xy}
\theoremstyle{plain}
\newtheorem{theorem}{Theorem}[section]
\newtheorem{lemma}[theorem]{Lemma}

\theoremstyle{definition}
\newtheorem{definition}[theorem]{Definition}
\newtheorem{notation}[theorem]{Notation}
\theoremstyle{remark}
\newtheorem*{remark}{Remark}

\newcommand{\C}{\ensuremath{\mathfrak{C}}}

\newcommand{\N}{\ensuremath{\mathfrak{N}_\mathcal{D}}}
\newcommand{\Real}{\ensuremath{\mathbb{R}}}

\newcommand{\Set}{\ensuremath{\mathbf{Set}}}
\newcommand{\ang}[1]{\ensuremath{\langle#1\rangle}}

\newcommand{\set}[1]{\ensuremath{\left\{#1\right\}}}
\newcommand{\such}[1]{\ensuremath{\left|#1\right.}}
\newcommand{\ds}[1]{\ensuremath{_{_{#1}}}}
\begin{document}
  \title[The continuity of the map lim.]{The continuity of the map $\lim_{T}$ in Hausdorff spaces}
  \author{J. E. Palomar Taranc\'on}
  \address{Dep. Math. Inst. Jaume I\\C./ Europa, 3-2E\\112530-Burriana-(Castell\'on)-Spain}
  \email{jpalomar016n@cv.gva.es}
 \secondemail{jepalomar23@ono.com}

\begin{abstract}
  Consider a Hausdorff space (X,T) and a set C of converging nets in
  X. By virtue of the limit uniqueness, the relation Lim which assigns
  each member x of X to every net N lying in C that converges to x is
  a map.  Of course, structuring C with some topology U, Lim can be a
  continuous map. If T is a topology induced by a uniformity, and F is
  a function space such that X is the codomain of each f in F, it is a
  well-known property the uniform limits of continuous functions to be
  continuous. In this paper, the author shows that the continuity of
  limits of continuous-map nets is implied by the continuity of the
  map Lim, whenever the involved topologies are large enough, being
  this result obtained without using neither the uniform convergence
  notion nor the uniformity concept.
\end{abstract}

\primaryclass{54C05} \secondaryclass{54B30, 54C35}
\keywords{Continuity, limits, concrete categories, partial-morphisms,
  partial algebras, net-functor.}  \maketitlepage
 \section[Introduction]{Introduction}
 \label{sec:introduction}
 Consider any non empty set $X$ and a topology $T$ for $X$. If
 $(X,T)$ is a Hausdorff space and $\C(T)$ stands for a collection of
  converging nets in $X$, then the relation $\lim_T\subseteq\C(T)\times{}X$ such that
  $(S,x)\in\lim_T$ implies $\lim_T S=x$ is a map. It is a natural question to
 consider some topology $\mathfrak{T}$ for $\C(T)$ such that
 \mbox{$\lim_T:\left(\C(T),\mathfrak{T}\right)\rightarrow{}(X,T)$} is a continuous
 map. Another natural question consists of knowing under what
 conditions  the limit $f$ of a net of continuous maps $S=(f_\delta)_{\delta\in D}$
  is again a continuous one. The second question possesses a
 well-known answer for uniform spaces, namely, the uniform convergence. In this paper we shall
 prove (Theorem~\ref{thm:main}) that, under some circumstances, the first property implies
 the second one.
 This result has been performed without using neither the uniform
 convergence notion  nor the uniformity concept.
 
 Of course the continuity of a map $f:(X_1,T_1)\rightarrow{}(X_2,T_2)$
 can be defined by means of an expression of the form $\lim_{T_2}
 f(x_\rho)= f(\lim_{T_1} x_\rho)$, or simply $\lim_{T_2} f=f
 \lim_{T_1}$; accordingly, if $(f_\delta)_{\delta\in D}$ is a net of
 continuous maps converging to $f$, then the continuity of the limit
 $f$ can be stated by the relation $\lim_{T_2} f= f \lim_{T_1}$,
 explicitly $\lim_{T_2} \lim_\delta f_\delta= \lim_\delta f_\delta
 \lim_{T_1}$.  Thus, since each of the $f_\delta$ is supposed to be
 continuous, the last equation implies that $\lim_{T_2} \lim_\delta
 f_\delta= \lim_\delta \lim_{T_2}f_\delta$, and this equation can be
 obtained from $\lim_{T_2} \lim_\delta = \lim_\delta \lim_{T_2}$
 whenever the last one is true.  Notice, that if $\lim_{T_2}$ is a
 map, then the last relation denotes $\lim_{T_2}$ to be continuous.
 This is why, at least in Hausdorff spaces, in order to state the
 continuity of the limits of converging continuous map nets, the
 continuity of the map $\lim_{T_2}$ is the natural condition, and such
 a requirement can be stated without using the uniform space concept.
 
 Finally, in Theorem~\ref{diff} it is shown that the continuity of the
 map $\lim_{T}$ implies some differential operators to be also
 continuous.

\section[The net-functor]{The net-functor and related partial algebras}

Let $\mathcal{D}$ be a small-class (set) of directed sets being stable
under products, that is, $\mathcal{D}$ contains with every countable
subset $\set{(D_n,\leq_n)\in\mathcal{D}\such n\in I\subseteq N}$ the
product directed set $\prod_{n\in I}(D_n,\leq_n)$, besides, for every
$(D,\leq)\in\mathcal{D}$ and each map $f:D \rightarrow{}\mathcal{D}$
the product directed set $\prod_{d\in D}f(d)$ again belongs to
$\mathcal{D}$. With these assumptions, denote by $\N:\Set
\rightarrow{}\Set$ the endofunctor in the category $\Set$ of ordinary
sets, the object-map ${\N}_O$ of which carries each set $X$ into the
net set $\bigcup_{D\in\mathcal{D}}X^D$; and the arrow-map ${\N}_A$
carries each morphism $f:X_1 \rightarrow{}X_2$ into the map $\N(f)(S)$
such that $(x_\delta)_{\delta \in
  D}\mapsto\left(f(x_\delta)\right)_{\delta\in D}\in{}X_2^D$.
\nocite{Adamek}

\begin{remark}
  Henceforth, to improve the readability, whenever a morphism $f$
  occurs as the argument of a functor $\N$, and in the same expression
  occurs an element $x$ as the argument for $\N(f)$, we shall write
  $\N\ang{f}(x)$ instead of $\N{}f(x)$ or $\N(f)(x)$, that is, the
  symbol $\ang{\,}$ in the expression $\N\ang{\,}$ indicates that $\N$
  works as an arrow-map.
\end{remark}

\begin{definition}
  Given a partial morphism \nocite{Rosolini} $\N(X)
  \xleftarrow{i_{\C(T)}}\C(T) \xrightarrow{\lambda_T}X$, where
  $i_{\C(T)}$ stands for the canonical inclusion and $\C(T)$ is any
  non-empty subset of $\N(X)$, say the pair
  $\left(X,(i_{\C(T)},\lambda_T)\right)$ to be a $\lambda$-partial
  $\N$-algebra, provided that the following condition is satisfied.

\begin{it}
  The set $\set{(S,x)\in\C(T)\times{}X\such x=\lambda(S)}$ is a
  convergence-class such that the corresponding topology is separated.
\end{it}

Obviously, in the category $\Set$ of ordinary sets and maps, the
canonical inclusion of a subset is a monomorphism, therefore so is
$i_{\C(T)}:\C(T) \rightarrow{}\N(X)$, and for every map
$\lambda:\C(T)\rightarrow{}Y$, the pair $(i_{\C(T)},\lambda_T)$ is a
partial map or partial morphism between sets.
\end{definition}

\begin{lemma}\label{lem:lim}
If $\left(X,(i_{\C(T)},\lambda_T)\right)$ is a $\lambda$-partial $\N$-algebra, then for
every $S$ in $\C(T)$ the relation $\lambda(S)=\lim_T S$ holds.
\end{lemma}

\begin{proof}
It is a straightforward consequence of the former definition.
\end{proof}

\begin{remark}
  Although the considered directed set class $\mathcal{D}$ is a small
  one, that is, a set of directed sets, in general, it is sufficient.
  Convergence classes can be build over those directed set collections
  each member of which is isomorphic to a base for a neighbourhood
  system ordered by inclusion. For instance, if a topological space
  satisfies the first axiom of countability only countable directed
  sets are necessary.
\end{remark}

From now on, a pair of the form $\left(X,(i_{\C(T)},\lim_T)\right)$
will be interpreted implicitly as a $\lambda$-partial $\N$-algebra,
therefore it must be understood $(X,T)$ to be a Hausdorff space,
otherwise $\lim_T$ need not be a map.

As usual, morphisms between algebras induced by endofunctors in $\Set$
are those maps with respect to which some diagrams commute. Thus,
given two $\lambda$-partial $\N$-algebras
$\left(X_2,(i_{\C(T_1)},\lim_{T_1})\right)$ and
$\left(X_2,(i_{\C(T_2)},\lim_{T_2})\right)$, a map
$f:X_1\rightarrow{}X_2$ is a morphism, provided that there is a unique
map $f^*$ such that the following diagram commutes.
\begin{equation}\label{eq:alg1}
\xymatrix{
&\ar[dl]_{i\ds{\C(T_1)}}\C(T_1)\ar[ddd]_{f^*}\ar[dr]^{\lim\ds{T_1}}&\\
\N(X_1)\ar[d]_{\N(f)}&&\ar[d]^{f}X_1\\
\N(X_2)&&X_2\\
&\ar[ul]^{i\ds{\C(T_2)}}\C(T_2)\ar[ur]_{\lim\ds{T_2}}& }
\end{equation}
Of course, the map $f^*$ is nothing but the restriction
$\N(f)\vert\ds{\C(T_1)}$ of $\N(f)$ to $\C(T_1)$; accordingly the
former diagram can be also written as follows.
\begin{equation}\label{eq:alg2}
\xymatrix{
&&\ar[dll]_{i\ds{\C(T_1)}}\C(T_1)\ar[ddd]_{\N(f)\vert\ds{\C(T_1)}}\ar[dr]^{\lim\ds{T_1}}&\\
\N(X_1)\ar[d]_{\N(f)}&&&\ar[d]^{f}X_1\\
\N(X_2)&&&X_2\\
&&\ar[ull]^{i\ds{\C(T_2)}}\C(T_2)\ar[ur]_{\lim\ds{T_2}}& }
\end{equation}

\begin{lemma}\label{lem:continuous}
Morphisms among $\lambda$-partial $\N$-algebras are, precisely,  continuous maps.
\end{lemma}

\begin{proof}
Taking into account (\ref{eq:alg2}), every morphism satisfies the relation
\nocite{Kelley}
\begin{equation}\label{eq:morph1}
  \lim_{T_2}\N\ang{f}(S)=f(\lim_{T_1}S)
\end{equation}
for every converging net $S=(x_\delta)_{\delta\in D}$ in $\C(T_1)$, and this expression
being written in the usual notation is nothing but
$\lim\ds{T_2}f(x_\delta)=f(\lim\ds{T_1}x_\delta)$.
\end{proof}

\begin{remark}
For our purposes, it is a better notation $\lim_{T_2}\N\ang{f}(S)=f(\lim_{T_1}\!S)$ than
$\lim\ds{T_2}f(x_\delta)=f(\lim\ds{T_1}x_\delta)$, since in the first expression it is
denoted explicitly $\N\ang{f}$ to have a net as argument.
\end{remark}

Since in a Hausdorff space $(X,T)$, $\lim_T:\C(T)\rightarrow{}X$ is a
map, it is a natural question to build some separated topology
$\mathfrak{T}$ for $\C(T)$ with respect to which $\lim_T$ is a
continuous map. According to (\ref{eq:alg2}), the continuity for
$\lim\ds{T}$ can be stated by means of the following diagram.
\begin{equation}\label{eq:alg3}
\xymatrix{
&&\ar[dll]_{i\ds{\C(\mathfrak{T})}}\C(\mathfrak{T})\ar[ddd]_{\N(\lim_T)\vert\ds{\C(\mathfrak{T})}}\ar[dr]^{\lim\ds{\mathfrak{T}}}&\\
\N\left(\C(\mathfrak{T})\right)\ar[d]_{\N(\lim_T)}&&&\ar[d]^{\lim\ds{T}}\C(T)\\
\N\left(X)\right)&&&X\\
&&\ar[ull]^{i\ds{\C(T)}}\C(T)\ar[ur]_{\lim\ds{T}}& }
\end{equation}
 Thus, $\lim_T:\C(T)\rightarrow{}X$ is a continuous map if and only if the former diagram commutes,
 that is, the following relation holds.
\begin{equation}\label{eq:limcont}
 \lim_T\left(
 \lim_{\mathfrak{T}}\right)=\lim_T\left(\N\left(\lim_T\right)\right)
\end{equation}

\subsection{Transposing a net of nets}

When the underlying set of a $\lambda$-partial $\N$-algebra
$\left(\C(T),(i_{\C(T)},\lim_\mathfrak{T})\right)$ is a net set
$\C(T)$ with a separate topology $\mathfrak{T}$, each net lying in
$\C(\mathfrak{T})$ is a net of nets. For example, given a directed set
$(D,\leq)\in\mathcal{D}$, let $S=(S_\delta)_{\delta\in D}$ be a net
lying in $\C(\mathfrak{T})$. Of course, for each $\delta\in{}D$,
$S_\delta=(x_{d,\delta})_{d\in D_\delta}$ is a net in $\C(T)$, where
$(D_\delta, \leq_\delta)$ stands for a directed set in $\mathcal{D}$,
for each $\delta\in D$.  Obviously, $S$ can be regarded as a net of
nets. If for every $\delta\in D$, the relation
\mbox{$(D_\delta,\leq_\delta)$} $=(D_0,\leq_0)$ holds, then the net
$S$ can be written in a matrix form, that is,
$S=(x_{d,\delta})_{(d,\delta)\in D_0\times D}$, being the nets of $S$
the $S_\delta$. However, transposing the matrix we obtain the net of
nets $S^t=(x_{\delta,d})_{(\delta,d)\in D\times D_0}$, which is formed
by the nets $S^t_d=(x_{\delta,d})_{\delta\in D}$, for each $d$ in
$D_0$.

\begin{definition}\label{defn:transposable} If a net of nets $S=(S_\delta)_{\delta\in D}$ belongs to
  $\C(\mathfrak{T})$ for some $\lambda$-partial $\N$-algebra
  $\left(\N(X),(i_{\C(\mathfrak{T})},\lim_{\mathfrak{T}})\right)$,
  then $\N\ang{\lim_{T}}(S)$ is the net $\displaystyle(x_\delta
  )_{\delta\in D}= (\lim_dx_{d,\delta})_{\delta\in D}=(\lim_T
  S_\delta)_{\delta\in D}$. Sometimes, the nets forming its transposed
  are also converging ones. In this case, there is also the net
  $\N\ang{\lim_T}(S^{t})$ which is formed by the nets $\displaystyle
  (x_d )_{d\in D_0}= (\lim_\delta x_{\delta,d})_{d\in
    D_0}=\N\ang{\lim_T} S^t$.  From now on, say a topology
  $\mathfrak{T}$ for a net set $\C(T)$ to satisfy the transposition
  property, provided that for every $\mathfrak{T}$-converging net of
  the form $S=(x_{d,\delta})_{d\in D_0,\delta\in D}$ there exists
  $\N\ang{\lim_T}S^t=(x_d)_{d\in D_0}$, besides, the following
  relation holds.
\begin{equation}\label{transpose}
\N\ang{\lim_T}S^t=\lim_\mathfrak{T}S
\end{equation}
\end{definition}
\begin{remark}
If $S$ is a one--column or one--row net, for instance
$S=(x_\delta)_{\delta\in D}$, both limits $\lim_T S$ and $\lim_T
S^t$ are the same; consequently one can eliminate the transposition
operator$(\,)^t$.
\end{remark}
\begin{remark}
The continuity of the map $\lim_T$, according to (\ref{eq:limcont})
is implied by the relation
\begin{equation}\label{eq:trans2}
  \lim_T \N\ang{\lim_T}S^t=\lim_T\lim_\mathfrak{T}S^t
\end{equation}
However, if we consider the net of nets $S$ in a matrix form, that
is, considering column-nets and row-nets, the former expression must
be written as follows
\begin{equation}\label{eq:trans3}
 \lim_T \N\ang{\lim_T}S^t=\lim_T\lim_\mathfrak{T}S
\end{equation}
because the net consisting of the limits of all rows is a column
one, and vice versa.
\end{remark}

\subsection{Continuity of the map $\lim_T$}
The results in this section are implied by the continuity of the map
$\lim_T$ together with the transposition property of the involved
topologies.
\begin{theorem}\label{trans}
  If $(X_1,(i_{\C(T_1)},\lim_{T_1}))$ and
  $(X_1,(i_{\C(T_1)},\lim_{T_1}))$ are two $\lambda$-partial algebras,
  and $\mathfrak{T}_1$ and $ \mathfrak{T}_2$ two topologies for
  $\C(T_1)$ and $\C(T_2)$, respectively, each of which satisfies the
  transposition property, then for every continuous map $f:(X_1,T_1)
  \rightarrow{}(X_2,T_2)$ the image $\N(f):\N(X_1,T_1)
  \rightarrow{}\N(X_2,T_2)$ is again continuous.
\end{theorem}
\begin{proof}
  Let $S=(S_d)_{d\in D_1}\in\N\left(\N(X)\right)$ be a converging net,
  such that for every $d\in D_1$, the net
  $S_d=(x_{\delta,d})_{\delta\in D_2}\in\N(X)$ is also a converging
  one. Since $\lim_{T_2}$ is a map, therefore a $\Set$-morphism, by
  the morphism composition compatibility of any functor it follows,
  that
\begin{equation}\label{thm:300}
\N\ang{\lim_{T_2}}\left(\N\ang{\N\ang{f}}(S)\right)^t=\N\left(\lim_{T_2}\left(\N\ang{f}\right)(S)^t\right)
\end{equation}

Now, by virtue of the transposition property the net of nets
$\left(\N\ang{\N\ang{f}}(S)\right)$ satisfies (\ref{transpose}), consequently

\begin{multline}\label{thm:30}
 \N\ang{\lim_{T_2}}\left(\N\ang{\N\ang{f}}(S)\right)^t=\\
 \N\ang{\lim_{T_2}}\N\ang{\N\ang{f}}(S)^t=\lim_{\mathfrak{T}_2}
\N\ang{\N\ang{f}}(S)
\end{multline}
hence, from (\ref{thm:300}) and (\ref{thm:30}) we obtain
\begin{equation}\label{thm:301}
\N\left(\lim_{T_2}\left(\N\ang{f}\right)(S)^t\right)=\lim_{\mathfrak{T}_2}\N\ang{\N\ang{f}}(S)
\end{equation}
and because $f$ is assumed to be continuous, taking into account
Lemma~\ref{lem:continuous} and equation (\ref{eq:morph1}) one can obtain that
$\displaystyle\N\left(\lim_{T_2}\left(\N\ang{f}\right)(S)^t\right)=\N\left(f(\lim_{T_1}S)^t\right)$,
and using this equality in (\ref{thm:301}),
\begin{equation}\label{thm:302}
\N\left(f(\lim_{T_1}S)^t\right)=\N\ang{f}\left(\N(\lim_{T_1})S^t\right)
=\lim_{\mathfrak{T}_2}\N\ang{\N\ang{f}}(S)
\end{equation}
Finally, using (\ref{transpose}) in the former equation, we have that
\begin{equation}\label{thm305}
   \N\ang{f}\left(\lim_{\mathfrak{T}_1}S\right)
=\lim_{\mathfrak{T}_2}\N\ang{\N\ang{f}}(S)
\end{equation}
 and by virtue of (\ref{eq:morph1}), the former equation implies
$\N\ang{f}=\N(f)$ to be continuous.
\end{proof}

\begin{notation}
Henceforth, to avoid any confusion in expressions containing iterated limits like the
following one
\[
\lim_{T_1}\lim_{T_2}(x_{\delta,d})\ds{(\delta,d)\in D_1\times D_2}
\]
we shall write the generic member $\delta$ of the corresponding directed set together
with the symbol denoting the topology, being both symbols separated by a semicolon. Thus,
the former expression will be written explicitly as follows.
\[
\lim_{T_1;\,d}\lim_{T_2;\,\delta}(x_{\delta,d})\ds{(\delta,d)\in D_1\times D_2}
\]
\end{notation}

\begin{lemma}\label{pointwise}
  Let $\left(X_1,(i_{\C(T_1)},\lim_{T_1})\right)$ and
  $\left(X_2,(i_{\C(T_2)},\lim_{T_2})\right)$ be two $\lambda$-partial
  $\N$-algebras and $T$ a topology for
  $\mathcal{F}=\hom\left((X_1,T_1),(X_2,T_2) \right)$. If $T$ is finer
  than or equivalent to the topology of the pointwise convergence,
  then the restriction of the arrow-map of $\N$ to $\mathcal{F}$ is
  continuous.
\end{lemma}
\begin{proof}
Let $S=(f_\delta)_{\delta\in D_1}$ a net in $\mathcal{F}$ converging to $f$, so then for
every net $N=(x_\rho)_{\rho \in D_2}$ we have that
\begin{equation}\label{thm:p1}
  \lim_T\N\ang{S}(N)=\lim_T \N\ang{f_\delta}_{\delta\in D_1}(N)=\lim_{T;\,\delta}
  \N\ang{f_\delta}(x_\rho)_{\rho\in D_2}
\end{equation}
and because $T$ is assumed to be finer than or equivalent to the topology of pointwise
convergence, for every $\rho\in D_2$ the relation
$\displaystyle\lim_{T;\,\delta}f_\delta(x_\rho)=\lim_{T_2;\,\delta}f_\delta(x_\rho)=f(x_\rho)$
holds, consequently
\[\displaystyle\lim_{T;\,\delta}\N\ang{f_\delta}(x_\rho)_{\rho\in
D_2}=\lim_{T_2;\,\delta}\N\ang{f_\delta}(x_\rho)_{\rho\in D_2}=\N\ang{f}(x_\rho)_{\rho\in
D_2}\] that is to say, $\displaystyle\lim_{T}\N\ang{S}(N)=\N\ang{f}(N)=\N\ang{\lim_T
S}(N)$, therefore
\begin{equation}\label{thm:p2}
\lim_T\N\ang{S}=\N\ang{\lim_T S}
\end{equation}
\end{proof}

\begin{notation}
For every couple of  sets $X$ and $Y$, denote by $ev:Y^{X}\times{}X\rightarrow{}Y$ the
general evaluation map, that is to say, $\forall(f,x)\in Y^X\times{}Y$, $ev(f,x)=f(x)$.
Thus, fixed a point $x\in{}X$, $ev_x(f)=ev(f,x)=f(x)$ is the ordinary evaluation map. By
symmetry, fixed a map $f$, denote  $ev_f:X\rightarrow{}Y$ by co-evaluation map,
although it is nothing but the same $f:X\rightarrow{}Y$.
\end{notation}

\begin{theorem}[Main theorem]\label{thm:main}
  Let $\left(X_1,(i_{\C(T_1)},\lim_{T_1})\right)$,  $\left(X_2,(i_{\C(T_2)},\lim_{T_2})\right)$ and
   $\left(F,(i_{\C(T)},\lim_{T})\right)$ be three $\lambda$-partial $\N$-algebras,
 and $\mathfrak{T}_1$, $\mathfrak{T}_2$ and $\mathfrak{T}$  topologies for
   $\C(T_1)$,   $\C(T_2)$ and $\C(T)$ respectively. If $F$ belongs to $X_2^{X_1}$,
    $F_0$ stands for any non-empty subset of $F$, and the following
  statements hold,
  \begin{description}
\item[a)] Each of the topologies $\mathfrak{T}_2$, $\mathfrak{T}_2$ and $\mathfrak{T}$
satisfies the transposition property.
\item[b)] $T$ is finer than or equivalent to the pointwise convergence
  topology.
\item[c)] Every map in $F_0$ is continuous, and so is every evaluation map $ev(f,x)$ for
every $x\in X_1$.
\item[d)] The map $\lim_{T_2}:\left(\C(T_2),\mathfrak{T}_2\right)\rightarrow{}(X_2,T_2)$ is continuous.
\end{description}
then, for every directed set $(D_1,\preceq)$ in $\mathcal{D}$, and each converging net
$S=(f_d)_{d\in D_1}$ in $F_0$ the limit
$f=\lim_\mathcal{T}S=\lim_T(f_\delta)\ds{\delta\in D_1)} $ is a continuous map.
\end{theorem}

\begin{proof}
Let  $N=(x_\rho)\ds{\rho\in D_2}$ a converging net in $X_1$ and consider the net of nets
\begin{equation}\label{eq:thm11}
S=\N\ang{\N\ang{ev}}\left(\N\ang{f_\delta}_{\delta\in{}D_1},(x_\rho)_{\rho\in{}D_0}\right)
\end{equation}
 Since, by hypothesis, the map $\lim_{T_2}$ is
continuous, from (\ref{eq:trans3}) it follows that
\begin{equation}\label{eq:thm2}
 \lim_{T_2}\N\ang{\lim_{T_2}}S^t=\lim_{T_2}\lim_{\mathfrak{T}_2}S
\end{equation}

Now, since it is assumed that, for each $f\in F_0$ the map $ev_f=f$ is
continuous, and by hypothesis, so is each $ev_x$ for every $x\in X_1$
together with their images under $\N$ and $\N\circ\N$, as consequence
of Theorem~\ref{trans}. Thus, by continuity, from the left hand of
(\ref{eq:thm2}) we obtain that

\begin{multline}\label{thm3}
 \lim_{T_2}\N\ang{\lim_{T_2}}S^t=\\
  \lim_{T_2;\,\delta}\,\N(\lim_{T_2;\,\rho})\left(\N\ang{\N\ang{ev}}\left(\left(\N\ang{f_\delta}\right)_{\delta\in D_1},(x_\rho)\ds{\rho\in D_2}\right)\right)^t\\=
  \lim_{T_2;\,\delta}\left(\N\ang{ev}\left((f_\delta)_{\delta\in D_1},\lim_{T_1;\,\rho}(x_\rho)_{\rho\in
  D_2}\right)\right)^t=
\lim_{T_2;\,\delta}\left(f_\delta\left(\lim_{T_1}x_\rho\right)\right)^t
 \end{multline}
and because, by assumption, $T$ is finer than or equivalent to the pointwise convergence
topology, the net $(f_\delta)_{\delta\in D_1}$ pointwise converges to $f$, therefore
\begin{equation}\label{thm3b}
  \lim_{T_2;\,\delta}\left(f_\delta\left(\lim_{T_1}x_\rho\right)\right)^t =
  \lim_{T;\,\delta}\left(f_\delta\left(\lim_{T_1}x_\rho\right)\right)^t =f(\lim_{T_1}x_\rho)
\end{equation}
where the transposition is eliminated because of the expression
contains only a one--row net.

Likewise, taking into account that $\mathfrak{T}_2$ satisfies the transposition property,
using (\ref{transpose}) in (\ref{eq:thm2}), from the right hand of (\ref{eq:thm2}) we
obtain that
\begin{multline}\label{thm:add}
\lim_{T_2}\lim_{\mathfrak{T}_2}S=
\lim_{T_2;\,\rho}\,\,\lim_{\mathfrak{T}_2;\,\delta}\N\ang{\N\ang{ev}}\left(\left(\N\ang{f_\delta}\right)_{\delta\in
D_1},(x_\rho)\ds{\rho\in
 D_2}\right)=\\
  \lim_{T_2;\,\rho}\N\ang{\lim_{T_2;\,\delta}}\left(\N\ang{\N\ang{ev}}\left(\N\ang{f_\delta}_{\delta\in D_1},(x_\rho)\ds{\rho\in
  D_2}\right)\right)^{t}
\end{multline}
Now, by continuity it follows that
\begin{multline}\label{eq:commuta2}
 \lim_{T_2;\,\rho}\N\ang{\lim_{T_2;\,\delta}}\left(\N\ang{\N\ang{ev}}\left(\N\ang{f_\delta}_{\delta\in D_1},(x_\rho)\ds{\rho\in
  D_2}\right)\right)^t=\\
\lim_{T_2;\,\rho}\left(\N\ang{ev}\left(\N\ang{\lim_{T_2;\,\delta}}\N\ang{f_\delta}_{\delta\in
D_1},(x_\rho)\ds{\rho\in
  D_2}\right)\right)^t=\\
  \lim_{T_2;\,\rho}\left(\N\ang{ev}\left(\N\ang{\lim_{T_2;\,\delta}f_\delta}_{\delta\in
D_1},(x_\rho)\ds{\rho\in
  D_2}\right)\right)^t
\end{multline}
and because $T$ is assumed to be finer than or equivalent to the topology of pointwise
convergence,
\begin{multline}\label{eq:commuta3}
\lim_{T_2;\,\rho}\left(\N\ang{ev}\left(\N\ang{\lim_{T_2;\,\delta}f_\delta}_{\delta\in
D_1},(x_\rho)\ds{\rho\in
  D_2}\right)\right)^t=\\
  \lim_{T_2;\,\rho}\left(\N\ang{ev}\left(\N\ang{\lim_{T;\,\delta}f_\delta}_{\delta\in
D_1},(x_\rho)\ds{\rho\in
  D_2}\right)\right)^t=\\\lim_{T_2;\,\rho}\left(\N\ang{ev}\left(\N\ang{f},(x_\rho)\ds{\rho\in
  D_2}\right)\right)^t=\\
\lim_{T_2;\,\rho}\left(\N\ang{f}\left((x_\rho)\ds{\rho\in D_2}\right)\right)^t=
\lim_{T_2;\,\rho}\N\ang{f}\left((x_\rho)\ds{\rho\in D_2}\right)
\end{multline}
 Finally, by virtue of equation (\ref{eq:thm2}), from (\ref{thm3b}) and
(\ref{eq:commuta3}) it follows immediately,
\begin{equation}\label{eq:thm5}
  \lim_{T_2}\N(f)\Big((x_\rho)_{\rho\in D_2}\Big)=f(\lim_{T_1}x_\rho)
\end{equation}
and because the net $(x_\rho)_{\rho\in D_2}$ is arbitrary,  $f$ is a continuous map.
\end{proof}

\subsection{The continuity of differential operators}
Let $\mathfrak{V}_1=\left(F_1,(i_{\C(T_1)},\lim_{T_1})\right)$ be a
$\lambda$-partial $\N$-algebra each member $\mathbf{f}$ of which is a
differentiable being its domain dense in a fixed closed subset $U$ of
$\Real$ with non-empty interior, and the codomain $\Real^n$.  Assume
$(F_1,T_1)$ to be a Hausdorff space. Of course, although the domains
of members of $F_1$ need not be the same, $(F_1,T_1)$ can be a
topological vector space, for there are separated topologies and
algebraic structures for such a kind of function sets, for instance,
see \cite{Palomar2}.  Let
$\mathfrak{V}_2=\left(F_2,(i_{\C(T_2)},\lim_{T_2})\right)$ another
$\lambda$-partial $\N$-algebra such that the associated topological
vector space $(F_2,T_2)$ contains the derivative of every member of
$(F_1,T_1)$. Obviously, the differential operator
$\dfrac{d}{dt}:(F_1,T_1) \rightarrow{}(F_2,T_2)$ is continuous,
provided that for each directed set $(D,\leq)$ in $\mathcal{D}$ and
every converging net $S=(\mathbf{f}_\delta)_{\delta\in D}$ the
following relation holds.
\begin{equation}\label{eq:diff}
  \lim_{T_2} \N\ang{\frac{d}{dt}}(S)=\frac{d}{dt}\lim_{T_1}S
\end{equation}

\begin{theorem}\label{diff}
  Let $\left(X\subseteq\Real,(i_{\C(T)},\lim_T)\right)$ and
  $\left(\Real^n,(i_{\C(T_u)},\lim_{T_u})\right)$ two
  $\lambda$-partial $\N$-algebras, where $T$ and $T_u$ are the
  standard topologies inducing the ordinary smooth structure for
  $\Real$ and $\Real^n$ respectively. Let $F_1\subseteq{}F_2$ be two
  subsets of $(\Real^n)^X$ such that each member of $F_1$ is defined
  and differentiable in a dense subset of $X$ and its derivative
  belongs to $F_2$. Let $\mathcal{T}$ stand for a separated topology
  for $F_2$ and $\mathcal{T}_0$ the relative one for $F_1$. Let
  $\left(F_1,(i_{\C(\mathcal{T}_0)},\lim_{\mathcal{T}_0})\right)$ and
  $\left(F_2,(i_{\C(\mathcal{T})},\lim_\mathcal{T})\right)$ two
  $\lambda$-partial $\N$-algebras. If the following statements hold,
\begin{description}
   \item[a)] $\mathcal{T}$ is finer than or equivalent to the
  pointwise convergence topology.
  \item[b)] The map $\lim_{T_u}$ is continuous.
\end{description}
then  $\dfrac{d}{dt}:(F_1,T_1) \rightarrow{}(F_2,T_2)$ is a continuous map.
\end{theorem}
\begin{proof}
  Consider a directed set $(D,\leq)$ in $\mathcal{D}$ and any
  converging net $S=(\mathbf{f}_\delta)_\delta\in D$ in $F_1$. By
  assumption every $\mathbf{f}_\delta$ is differentiable for any
  $\delta \in D$. First, for each fixed $p\in X$, define the map
  $\mathbf{g}_{\delta,p}(x)=\dfrac{\mathbf{f}_\delta(p+x)-\mathbf{f}_\delta{}(p)}{x}$.
  With these assumptions, if $N$ is any net in $X$ converging to $0$,
  then for every $p\in x$, taking into account Lemma~\ref{pointwise},
\begin{equation}\label{thm:diff1}
  \frac{d}{dt}\lim_{\mathcal{T}_0}(\mathbf{f}_\delta)_{\delta\in D}|_p=
  \lim_{T_u}\N\ang{\lim_{\mathcal{T}_0}\mathbf{g}_{\delta,p}}(N)=
  \lim_{T_u}\lim_{\mathcal{T}_0}(\N\ang{\mathbf{g}}_{\delta,p})_{\delta\in{}D}(N)
\end{equation}
Since $\lim_{T_u}$ by assumption is continuous, and taking into account statement
\textbf{a}),
\begin{multline}\label{thm:diff2}
\lim_{T_u}\lim_{\mathcal{T}_0}(\N\ang{\mathbf{g}}_{\delta,p})_{\delta\in{}D}(N)=
  \lim_{\mathcal{T}}\N\ang{\lim_{T_u}}(\N\ang{\mathbf{g}}_{\delta,p})_{\delta\in{}D}(N)=\\
  \lim_{\mathcal{T}}\left(\N\ang{\lim_{T_u}\mathbf{g}_{\delta,p}}\right)_{\delta\in{}D}(N)=
\lim_{\mathcal{T}}\N\ang{\frac{d}{dt}}(\mathbf{f}_\delta)_{\delta\in{}D}|_p
\end{multline}
and because $p$ is arbitrary, from (\ref{thm:diff1}) and (\ref{thm:diff2}) we obtain that
\begin{equation}\label{thm:diff3}
 \frac{d}{dt}\lim_{\mathcal{T}_0}(\mathbf{f}_\delta)_{\delta\in D}=
 \lim_{\mathcal{T}}\N\ang{\frac{d}{dt}}(\mathbf{f}_\delta)_{\delta\in{}D}
\end{equation}
\end{proof}


\begin{thebibliography}

\bibitem{Adamek}
\textbf{Jirí Ad\'{a}mek}, \textbf{Horst Herrlich}, \textbf{George~E Strecker},
  \emph{{Abstract and concrete categories. The joy of cats.}},
  {Wiley-Interscience Publication. New York etc.: John Wiley \& Sons, Inc. xii,
  482 p.} (1990)

\bibitem{Kelley}
\textbf{John~L Kelley}, \emph{General topology}, Springer-Verlag, New York
  (1975)

\bibitem{Palomar2}
\textbf{{J}~{E} Palomar~Taranc{\'o}n}, \emph{Vector spaces over function
  fields. {V}ector spaces over analytic function fields being associated to
  ordinary differential equations}, Southwest J. Pure Appl. Math.  (2000)
  60--87 (electronic)

\bibitem{Rosolini}
\textbf{Giuseppe Rosolini}, \emph{Domains and dominical categories}, Riv. Mat.
  Univ. Parma (4) 11 (1985) 387--397

\end{thebibliography}
\end{document}